\documentclass [11pt]{amsart}
\usepackage{amssymb}
\usepackage{xypic}

\newcommand{\dual}{\makebox[0mm]{}^{{\scriptstyle\vee}}}
\newcommand{\ddual}{\makebox[0mm]{}^{{\scriptstyle\vee\vee}}}

\newtheorem{theorem}{Theorem}[section]

\newtheorem{lemma}[theorem]{Lemma}
\newtheorem{proposition}[theorem]{Proposition}
\newtheorem{definition}[theorem]{Definition}
\newtheorem{corollary}[theorem]{Corollary}

\newtheorem{exmp}[theorem]{Example}
\newtheorem{exmps}[theorem]{Examples}
\newtheorem{rem}[theorem]{Remark}
\newenvironment{example}{\begin{exmp}\rm}{\end{exmp}}

\newenvironment{remark}{\begin{rem}\rm}{\end{rem}\rm}

\newcommand{\beeq}[1]{\begin{eqnarray}\label{#1}}
\newcommand{\eneq}{\end{eqnarray}}
\renewcommand{\to}{\xymatrix@1@=15pt{\ar[r]&}}
\renewcommand{\mapsto}{\xymatrix@1@=15pt{\ar@{|->}[r]&}}
\renewcommand{\twoheadrightarrow}{\xymatrix@1@=15pt{\ar@{->>}[r]&}}
\newcommand{\congpf}{\xymatrix@1@=15pt{\ar[r]^-\sim&}}
\renewcommand{\hookrightarrow}{\xymatrix@1@=15pt{\ar@{^{(}->}[r]&}}

\newcommand{\IC}{{\mathbb C}}

\newcommand{\IP}{{\mathbb P}}
\newcommand{\IQ}{{\mathbb Q}}
\newcommand{\IR}{{\mathbb R}}
\newcommand{\IZ}{{\mathbb Z}}

\newcommand{\Pic}{{\rm Pic}}
\newcommand{\Hom}{{\rm Hom}}
\newcommand{\Ext}{{\rm Ext}}

\newcommand{\Coh}{{\rm Coh}}

\newcommand{\coker}{{\rm coker}}
\newcommand{\ext}{{\rm ext}}
\newcommand{\End}{{\rm End}}

\newcommand{\Ker}{{\rm ker}}

\newcommand{\rk}{{\rm rk}}
\newcommand{\Db}{{\rm D}^{\rm b}}

\newcommand{\kc}{{\mathcal C}}

\newcommand{\ka}{{\mathcal A}}

\newcommand{\kh}{{\mathcal H}}

\newcommand{\ke}{{\mathcal E}}
\newcommand{\kf}{{\mathcal F}}

\newcommand{\kp}{{\mathcal P}}

\newcommand{\ko}{{\mathcal O}}
\newcommand{\kt}{{\mathcal T}}

\newcommand{\kx}{{\mathcal X}}

\newcommand{\verylongarrow}[1]{\hbox to #1{\rightarrowfill}}

\begin{document}

\title{Derived and abelian equivalence of K3 surfaces}
\author{Daniel Huybrechts}
\maketitle

Any isomorphism $X\cong X'$ between two schemes induces an
equivalence $\Coh(X)\cong\Coh(X')$ between their abelian
categories of coherent sheaves. Due to a classical result of
Gabriel \cite{Gabriel} the converse holds true as well. Thus,
$$\xymatrix{X\cong X'~~\ar@{<=>}[r]&~~\Coh(X)\cong\Coh(X').}$$

Is there a similar statement when isomorphism of schemes is
replaced by derived equivalence? More precisely, can one
naturally associate abelian categories to two varieties $X$ and
$X'$ such that the varieties are derived equivalent if and only if
there exists an equivalence between the abelian categories?

We will restrict to complex K3 surfaces and prove

\begin{theorem}\label{ThmfirstIntro}
  Two complex projective K3 surfaces $X$ and $X'$ are derived equivalent if
and only if there exist complexified K\"ahler classes $B+i\omega$
and $B'+i\omega'$ on $X$ respectively $X'$ such that the two
abelian categories $\ka_X(\exp(B+i\omega))$ and
$\ka_{X'}(\exp(B'+i\omega'))$ are equivalent. Thus,
$$\xymatrix{\Db(X)\cong\Db(X')~~\ar@{<=>}[r]&~~\ka_X(\exp(B+i\omega))\cong\ka_{X'}(\exp(B'+i\omega')).}$$
\end{theorem}

 Here, $\Db(X)$ is the bounded derived category $\Db(\Coh(X))$
and the equi\-va\-lence on the left hand side is linear and
exact. By definition  $\ka(\exp(B+i\omega))$ is the full
subcategory of all complexes $F^\bullet\in\Db(X)$ with cohomology
concentrated in degree $-1$ and zero and such that
$\kh^{-1}(F^\bullet)$ is torsion free with $\mu_{\rm
max}\leq(B.\omega)$ and the torsion free part of
$\kh^0(F^\bullet)$ satisfies $\mu_{\rm min}>(B.\omega)$. (For the
notation and the details of this definition see Section
\ref{SectAEyieldsDE}.) In our situation, $B+i\omega$ can be taken
as a complexified ample class, i.e.\
$B+i\omega\in\Pic(X)\otimes\IC$ with $\omega\in\Pic(X)$ ample. The
two directions of the theorem are proved in Section
\ref{SectAEyieldsDE} (Cor.\ \ref{DerEquiv}) respectively Section
\ref{SectDA} (Cor.\ \ref{CorDEyieldsAE}).

 The abelian category $\ka_X(\exp(B+i\omega))$ is the heart of a
$t$-structure on $\Db(X)$ that has been studied by Bridgeland in
\cite{BK3} in his quest for stability conditions on the
triangulated category $\Db(X)$. Bridgeland introduced the concept
of stability conditions, in \cite{B} in an effort to understand
Douglas's work on stability of branes. Roughly, a stability
condition on a triangulated category consists of a $t$-structure
and a stability function on its heart satisfying the
Harder--Narasimhan property. It is surprisingly difficult to
construct stability conditions on the derived category of a
higher dimensional projective Calabi--Yau variety. In \cite{BK3}
Bridgeland considers the case of K3 surfaces. In order to get
started he needs to construct explicit examples of stability
conditions and the heart of those are the abelian categories
$\ka(\exp(B+i\omega))$. (The abelian category $ \Coh(X)$ never
occurs as the heart of a stability condition.)

This paper grew out of the attempt to understand the geometric
meaning of the abelian categories $\ka(\exp(B+i\omega))$.

The naive idea to prove Gabriel's result that any equivalence
$\Coh(X)\cong\Coh(X')$ induces an isomorphism of the K3 surfaces
$X$ and $X'$ is the following. First note that the \emph{simple}
objects (or the \emph{minimal} objects, as we will call them, see
Section \ref{SectAEyieldsDE}) of $\Coh(X)$ are the structure
sheaves $k(x)$ of closed points $x\in X$. Since this is an
intrinsic notion, any equivalence of abelian categories sends
minimal objects  to minimal objects. Thus, the equivalence
$\Coh(X)\cong\Coh(X')$ induces a bijection $X\cong X'$ and, in
order to fully prove Gabriel's result, one only has to show that
this bijection is a morphism.

From this point of view it is natural to wonder how the minimal
objects of $\ka(\exp(B+i\omega))$ look like. This question is
interwoven with Theorem \ref{ThmfirstIntro} and we shall give
 the following complete classification in Section \ref{SectMinObj}
 (Prop.\ \ref{classmin}):

 \begin{theorem} For a K3 surface $X$
the minimal objects in $\ka(\exp(B+i\omega))$ are precisely the
objects
  \begin{itemize}
  \item[$\bullet$] $k(x)$, where $x\in X$ is a closed point and
  \item[$\bullet$] $F[1]$, where $F$ is a $\mu$-stable locally free sheaf
  with $\mu(F)=(B.\omega)$.\end{itemize}
 \end{theorem}

Thus, any equivalence between $\ka_{X'}$ and $\ka_{X}$ will
either induce an isomorphism $X'\cong X$ or will map closed
points in $X'$ to shifted $\mu$-stable vector bundles on $X$. In
order to combine both theorems, we have to prove a stronger
version of Orlov's well-known result saying that two K3 surfaces
are derived equivalent if and only if one is a moduli space of
stable sheaves on the other. In Proposition \ref{FMviamu} we
actually prove that in Orlov's result one can replace `stable' by
`$\mu$-stable' and `sheaves' by 'vector bundles'.

The abelian category $\ka(\exp(B+i\omega))$ plays a decisive role
in D-equivalence of K3 surfaces, but it also appears naturally
from a differential-geometric point of view. The minimal objects
of $\ka(\exp(i\omega))$, besides the point sheaves, are (shifted)
hyperholomorphic bundles, i.e.\ bundles that are holomorphic with
respect to all hyperk\"ahler rotations (with respect to $\omega$)
of the original complex structure. A short discussion of this
point of view is included in Section \ref{SectTw}.

The last section of this paper proves stability of Fourier--Mukai
transforms of certain $\mu$-stable vector bundles. The main
result not only yields stability in cases not covered by existing
result, but it gives, maybe more interestingly, a conceptual
explanation when and why stability of a Fourier--Mukai transform
of a $\mu$-stable vector bundle can be expected with Mukai vector
$v=(r,\ell,s)$. There we prove

\begin{theorem}
There exists a polarization $H'$ on $X'$ such that for any
  $\mu$-stable vector bundle $E$ on $X$ with $\mu(E)=-(\ell.H)/r$
  one has either
  \begin{itemize}
  \item[$\bullet$]  $\Phi(E)\cong k(y)[-2]$ if $[E\dual]\in M_H(v)$ or otherwise
  \item[$\bullet$] $\Phi(E)\cong F[-1]$ with $F$ a $\mu_{H'}$-stable vector bundle on
  $X'$.
  \end{itemize}
\end{theorem}

{\bf Acknowledgements.} I am grateful to Tom Bridgeland for
writing \cite{BK3}, answering my questions on  it,  and reminding
me of some of the many useful results in Mukai's paper \cite{Mu}.
Discussions with Alastair King were helpful to  get my ideas
straight.


\section{Abelian equivalence yields derived
equivalence}\label{SectAEyieldsDE}

A \emph{torsion pair} in an abelian category $\kc$ is a pair of
full subcategories $\kt,\kf\subset\kc$ such that $\Hom_\kc(T,F)=0$
for all $T\in\kt$, $F\in\kf$ and such that every object $E\in\kc$
fits into a short exact sequence
$$\xymatrix{0\ar[r]&T\ar[r]&E\ar[r]&F\ar[r]&0}$$
for some $T\in\kt$ and $F\in\kf$.

Following \cite{HRS} one associates to a given torsion pair
$(\kt,\kf)$ in $\kc$ a $t$-structure on the bounded derived
category $\Db(\kc)$ by setting $${\rm
D}^{\leq0}:=\{F^\bullet\in\Db(\kc)~|~H^i(F^\bullet)=0,~i>0;\
H^0(F^\bullet)\in\kt\}$$ and $${\rm
D}^{\geq0}:=\{F^\bullet\in\Db(\kc)~|~H^i(F^\bullet)=0,~i<-1;\
H^{-1}(F^\bullet)\in\kf\}.$$ Its heart, also called the
\emph{tilt} of $\kc$, is the abelian category
\begin{eqnarray*}
&\ka(\kt,\kf):={\rm D}^{\leq0}\cap {\rm
D}^{\geq0}\\
&=\{F^\bullet\in\Db(\kc)~|~H^i(F^\bullet)=0, ~ i\ne0,-1;\
H^0(F^\bullet)\in\kt;\ H^{-1}(F^\bullet)\in\kf\}.
\end{eqnarray*}

Thus, any object in $\ka(\kt,\kf)$ is isomorphic to a complex of
the form $$\xymatrix{F^{-1}\ar[r]^\varphi &F^0}$$ with
$\coker(\varphi)\in\kt$ and $\Ker(\varphi)\in\kf$. Note that in
particular $\kf[1]$ and $\kt$ are both naturally contained in
$\ka(\kt,\kf)$. Moreover, $(\kf[1],\kt)$ is a torsion pair in
$\ka(\kt,\kf)$ whose tilt is $\kc[1]$.

A torsion pair $(\kt,\kf)$ is called \emph{tilting} if every
object in $\kc$ is a subobject of an object in $\kt$. Similarly,
$(\kt,\kf)$ is \emph{cotilting} if every object in $\kc$ is a
quotient of an object in $\kf$. In the latter case, every object
in $\kc$ admits a resolution of length two by objects in $\kf$.
Indeed, any subobject of an object in $\kf$ is in $\kf$.

Suppose $(\kt,\kf)$ is a cotilting torsion pair. Then the natural
inclusion $\kf\subset\kc$ induces an exact equivalence
$\Db(\kf)\to\Db(\kc)$ (see \cite[Lemma 5.4.2]{BvdB}). Similar, if
$(\kt,\kf)$ is tilting, then $\Db(\kt)\cong\Db(\kc)$ is an
equivalence. Using that $(\kt,\kf)$ is a cotilting torsion pair
in $\kc$ if and only of $(\kf[1],\kt)$ is a tilting pair in
$\ka(\kt,\kf)$ (see \cite[Prop.\ I.3.2]{HRS}) one obtains for a
cotilting pair $(\kt,\kf)$ two exact equivalences
$$\Db(\kf)\cong\Db(\kc)~~{\rm
and}~~\Db(\kf)\cong\Db(\kf[1])\cong\Db(\ka(\kt,\kf)).$$ This
yields

\begin{proposition}
  For any cotilting torsion pair $(\kt,\kf)$ in
  an abelian category $\kc$ there  exists an exact equivalence
  $$\Db(\kc)\cong\Db(\ka(\kt,\kf)).$$
\end{proposition}
The result was first proved in \cite{HRS} under additional
assumptions (e.g.\ the existence of enough injectives in $\kc)$
and in the above form in \cite[Prop.\ 5.4.3]{BvdB}.

\bigskip

Let us now turn to the more concrete situation where the abelian
category $\kc$ is the category ${\Coh}(X)$ of coherent sheaves on
a smooth projective variety $X$ of dimension $n$.

We fix a polarization $H$ or, more generally, a K\"ahler class
$\omega$, so that degree $\deg(F)$ and slope $\mu(F)$ of a
coherent sheaf $F$ on $X$ can be defined as $\deg(F):=\int_X {\rm
c}_1(F).H^{n-1}$ (respectively $\int_X{\rm c}_1(F).\omega^{n-1}$)
and $\mu(F):=\deg(F)/\rk(F)$.

The \emph{Harder--Narasimhan-filtration} (HN-filtration for short)
of a coherent sheaf $F$ is the unique filtration
$$0\subset F_0\subset F_1\subset\ldots\subset F_n=F$$
that satisfies the following conditions: i) $F_0$ is the torsion
part of $F$, ii) The quotients $F_{i+1}/F_i$ are torsion free and
$\mu$-semistable for $i=0,\ldots,n-1$, and iii)
$\mu(F_1/F_0)>\ldots>\mu(F_n/F_{n-1})$.

Existence and uniqueness are easy to prove, see e.g.\ \cite[Thm.\
1.6.7]{HL}. One denotes $\mu_{\rm max}(F):=\mu(F_1/F_0)$ and
$\mu_{\rm min}(F):=\mu(F_n/F_{n-1})$ provided $F$ is not torsion.

For $\beta\in \IR$ one introduces the full subcategories
$$\kt(\beta),\kf(\beta)\subset\Coh(X).$$ By definition $\kt(\beta)$
is the category of all coherent sheaves $F$ with $\mu_{\rm
min}(F)>\beta$ or $F$ is a torsion sheaf and $\kf(\beta)$ is the
category of all torsion free coherent sheaves $F$ with $\mu_{\rm
max}(F)\leq\beta$ or $F\cong0$.

The following is an immediate consequence of the basic properties
of $\mu$-semistable sheaves and the existence of the
HN-filtration.

\begin{proposition}
  With the above notation $\kt(\beta),\kf(\beta)\subset\Coh(X)$ is a
  torsion pair.\qed
\end{proposition}

Torsion pairs of this form have been introduced by A.\ Schofield
and were later  studied for curves in \cite{P} and for K3 surfaces
in \cite{BK3}.

Let us denote the heart of the induced $t$-structure by
$\ka(\beta)$ (or, $\ka_X(\beta)$ if the dependence on $X$ needs
to be stressed). Thus, $\ka(\beta)$ is the following full
subcategory of $\Db(X)=\Db(\Coh(X))$:
$$\ka(\beta):=\{\xymatrix{F^{-1}\ar[r]^-\varphi &F^0}~|~\ker(\varphi)\in\kf(\beta),~\coker(\varphi)\in\kt(\beta)\}.$$
Note that, although not reflected by the notation, $\ka(\beta)$
depends on $\beta$ and on the chosen polarization (respectively
K\"ahler class).

\begin{corollary}\label{DerEquiv}
  Suppose $X$ and $X'$ are smooth projective varieties endowed
  with polarizations $H$ respectively $H'$ (or
  K\"ahler classes $\omega$ respectively $\omega'$). If for two
  real numbers $\beta,\beta'$ the abelian categories
  $\ka_X(\beta)$ and $\ka_{X'}(\beta')$ are equivalent, then $X$
  and $X'$ are derived equivalent. In other words,
  $$\xymatrix@=22pt{\ka_X(\beta)\cong\ka_{X'}(\beta')
  \ar@{=>}[r]&\Db(X)\cong\Db(X').}$$
\end{corollary}

\begin{remark}
  HN-filtrations not only exist with respect to $\mu$-semistability
  (see \cite[Ch.\ 1]{HL}). As only the formal properties were
  used, the above discussion goes through unchanged for e.g.\
  Gieseker stability for which the slope $\mu(F)$ is replaced by
  the Hilbert polynomial $\chi(F(n))$.
\end{remark}

In the following we shall be interested in the case of an
algebraic K3 surface $X$. The K3 surface $X$ will be endowed with
a K\"ahler class $\omega\in{\rm NS}(X)\otimes\IR\cong
H^{1,1}(X,\IZ)\otimes\IR$ and a B-field $B\in{\rm
NS}(X)\otimes\IR$. Then let $\beta:=(B.\omega)$.

Instead of considering $B$ and $\omega$ or the complexified
K\"ahler class $B+i\omega$, one more naturally uses
$\exp(B+i\omega)=1+(B+i\omega)+(B+i\omega)^2/2\in H^*(X,\IC)$,
which can be seen as a generalized Calabi--Yau structure on $X$.
This point of view fits nicely with the picture proposed by mirror
symmetry. For a discussion see \cite{HCY}.

Changing the notation of \cite{BK3} slightly we shall thus write
$$\kt:=\kt(\exp(B+i\omega)):=\kt(\beta),~~\kf:=\kf(\exp(B+i\omega)):=\kf(\beta)$$
and $$\ka(\exp(B+i\omega)):=\ka(\beta).$$

In particular Corollary \ref{DerEquiv} for two K3 surfaces $X$
and $X'$ reads:
$$\xymatrix{\ka_X(\exp(B+i\omega))\cong\ka_{X'}(\exp(B'+i\omega'))
\ar@{=>}[r]&\Db(X)\cong\Db(X').}$$

Using the Mukai vector $v(F)={\rm ch}(F)\sqrt{{\rm
td}(X)}=(r,\ell,s)$, Bridgeland introduces
$$Z(F):=\langle v(F),\exp(B+i\omega)\rangle$$ in order to
construct a stability function on the abelian category
$\ka(\exp(B+i\omega))$ (see \cite[Sect.\ 5]{BK3}). Here,
$\langle~~,~~\rangle$ is the Mukai pairing.

Clearly, $${\rm Im}\left\langle
v(F),\exp(B+i\omega)\right\rangle=(\ell.\omega)-r(B.\omega)
=(\ell.\omega)-r\beta.$$ Thus, for $r\ne0$, the slope $Z(F)$ is
contained in the upper half plane if and only if $\mu(F)>\beta$.
Also, for $F\in\kt$ one has ${\rm Im}(Z(F))\geq0$ and, similarly,
if $F\in\kf$, then ${\rm Im}(Z(F))\leq0$. Using $v(F[1])=-v(F)$,
this shows that ${\rm Im}(Z(F^\bullet))\geq0$ for all
$F^\bullet\in\ka(\exp(B+i\omega))$.

\begin{remark}\label{RemStabBrid}
  In \cite[Lemma 5.2]{BK3} it is shown that
  $Z(F^\bullet)\in\IR_{>0}\exp(i\pi\phi(F^\bullet))$ with the
  phase $\phi(F^\bullet)$ satisfying $0<\phi(F^\bullet)\leq1$
  holds for
  all $0\ne F^\bullet\in\ka(\exp(B+i\omega))$ if and only if
  $Z(F)\not\in\IR_{\leq0}$ for all spherical sheaves $F$. Note
  that the latter holds as soon as $(\omega.\omega)>2$.
\end{remark}

For later use we note that ${\rm Im}(Z(k(x)))=0$ for any closed
point $x\in X$ and ${\rm Im}(Z(F[1]))=0$ for any $\mu$-stable
vector bundle $F$ with $\mu(F)=\beta$. Under the assumption of
the remark, this is equivalent to $\phi(k(x))=\phi(F[1])=1$.

\section{Minimal objects in $\ka$}\label{SectMinObj}

The aim of this section is to classify minimal objects in
$\ka(\exp(B+i\omega))$ (modulo a technical result postponed to the
next section).

Recall that a non-trivial object $A$ in an abelian category $\ka$
is called \emph{minimal} if any surjection $A\twoheadrightarrow B$
with $B\ne0$ is an isomorphism. Equivalently, $A$ is minimal if
and only if every injection $0\ne C~\hookrightarrow A$ is an
isomorphism, i.e.\ $A$ has no proper subobjects. Usually, objects
of this type are called \emph{simple}, but `simple' for a sheaf
$F\in\Coh(X)$ has also a different meaning, i.e.\ that
$\End(F)=k$, so we rather use `minimal' instead.

Here are a few easy observations. Suppose $A$ is minimal and
$\varphi:A\to B$ is a morphism. Then either $\varphi=0$ or
$\varphi$ is injective. If in addition $B$ is minimal as well,
then either $\varphi=0$ or $\varphi$ is an isomorphism.

\begin{example}\label{exCoh}
As was mentioned earlier in the introduction, the minimal objects
in $\Coh(X)$ are the point sheaves $k(x)$ with $x\in X$ a closed
point. These minimal objects have the additional property that any
non-trivial sheaf $F\in \Coh(X)$ admits a surjection
$F\twoheadrightarrow k(x)$ for some $x\in X$.

Note however that they do not generate $\Coh(X)$. Recall that a
collection of objects in an abelian category generates the
category if every object admits a filtration whose quotients are
isomorphic to objects in the collection.
\end{example}

We shall be interested in the abelian category
$\ka:=\ka(\exp(B+i\omega))$ on a K3 surface $X$, which is by
definition a full subcategory of the derived category $\Db(X)$
obtained as a tilt of $\Coh(X)$ with respect to a torsion pair
$(\kt,\kf)$.

As this will be frequently used in the following discussion, we
recall the following standard fact (see \cite[p.\ 415]{KS}): Let
$\ka$ be the heart of a $t$-structure on a triangulated category
${\rm D}$. If
$\xymatrix@=14pt{0\ar[r]&A\ar[r]&B\ar[r]&C\ar[r]&0}$ is a short
exact sequence in $\ka$, then there exists a map $C\to A[1]$ such
that $\xymatrix@=14pt{A\ar[r]&B\ar[r]&C\ar[r]&A[1]}$ is a
distinguished triangle in ${\rm D}$. Conversely, if $A\to B\to
C\to A[1]$ is a distinguished triangle in ${\rm D}$ with objects
$A,B,C$ in $\ka$, then $\xymatrix@=14pt{0\ar[r]& A\ar[r]& B\ar[r]&
C\ar[r]&0}$ is a short exact sequence in $\ka$.

\medskip

In the following $B+i\omega\in{\rm NS}(X)_\IC$ is a complexified
K\"ahler class, i.e.\ $\omega\in{\rm NS}(X)_\IR$ is a K\"ahler
class and $B\in{\rm NS}(X)_\IR$ is arbitrary. Stability is
considered with respect to $\omega$ and we do not assume $\omega$
or $B$ to be rational.

\begin{proposition}\label{classmin}
  The minimal objects in $\ka(\exp(B+i\omega))$ are precisely the objects
  \begin{itemize}
  \item[$\bullet$] $k(x)$, where $x\in X$ is a closed point and
  \item[$\bullet$] $F[1]$, where $F$ is a $\mu$-stable locally free sheaf
  with $\mu(F)=(B.\omega)$.\end{itemize}
\end{proposition}

\begin{proof} To shorten the notation we write $(\kt,\kf)$ for the
  torsion pair induced by $\exp(B+i\omega)$. Similarly,
  $\ka:=\ka(\exp(B+i\omega))$. As before, we use
  $\beta:=(B.\omega)$.

\smallskip

  {\bf i) Point sheaves are minimal.}
  First, $k(x)\in\kt\subset\ka$ for any closed point $x\in X$.
  Suppose $k(x)\twoheadrightarrow F^\bullet$ is a non-trivial
  surjection in $\ka$ which we complete to a short exact sequence
  $$\xymatrix{0\ar[r]&E^\bullet\ar[r]&k(x)\ar[r]&F^\bullet\ar[r]&0}$$
  in $\ka$. Considered as a distinguished triangle in $\Db(X)$ it
  yields the long exact sequence
  $$\xymatrix@=15pt{0\ar[r]&\kh^{-1}(E^\bullet)\ar[r]&0\ar[r]&\kh^{-1}(F^\bullet)\ar[r]&
  \kh^0(E^\bullet)\ar[r]&k(x)\ar[r]&\kh^0(F^\bullet)\ar[r]&0.}$$
  Thus, $E^\bullet\cong E$ with $E\in\kt$. Moreover,
  $\kh^{-1}(F^\bullet)$ and $\kh^0(E^\bullet)\cong E$ are
  isomorphic on $X\setminus\{x\}$. Since $\kh^{-1}(F^\bullet)\in\kf$,
  this yields for $\kh^{-1}(F^\bullet)\ne0$ the contradiction
  $$\beta<\mu_{\rm min}(E)=\mu_{\rm
  min}(\kh^{-1}(F^\bullet))\leq\mu_{\rm
  max}(\kh^{-1}(F^\bullet))\leq\beta.$$ Hence, $\kh^{-1}(F) =0$,
  i.e.\ $F^\bullet\cong\kh^0(F^\bullet)\ne0$. By minimality of
  $k(x)$ as an object in $\Coh(X)$, the surjection $k(x)\twoheadrightarrow
  \kh^0(F^\bullet)$ is an isomorphism and hence the surjection
  $k(x)\twoheadrightarrow F^\bullet$ in $\ka$ is one.

\smallskip

  {\bf ii) Stable vector bundles of slope $\beta$ are minimal.}
  Let $F$ be a $\mu$-stable locally free sheaf with
  $\mu(F)=\beta$. Then $F[1]\in\ka$ by definition of $\ka$.
  Consider a short exact sequence
  \begin{equation}\label{minses}
  \xymatrix{0\ar[r]&G^\bullet\ar[r]&F^\bullet[1]\ar[r]&E^\bullet\ar[r]&0}
  \end{equation}
  in $\ka$. In order to show that $F[1]\in\ka$ is minimal, one
  proves that either $G^\bullet=0$ or $E^\bullet=0$.
  The long exact cohomology sequence of (\ref{minses}) considered
  as a distinguished triangle in $\Db(X)$ reads
  $$\xymatrix@=15pt{0\ar[r]&\kh^{-1}(G^\bullet)\ar[r]&F\ar[r]^-\varphi&\kh^{-1}(E^\bullet)\ar[r]&\kh^0(G^\bullet)
  \ar[r]&0\ar[r]&\kh^0(E^\bullet)\ar[r]&0.}$$
  Hence, $E^\bullet\cong E[1]$, where $E$ is torsion free with
  $\mu_{\rm max}(E)\leq\beta$. Consider the morphism
  $\varphi:F\to\kh^{-1}(E^\bullet)\cong E$ and its image $E'$. If
  $\varphi$ is neither trivial nor injective, then $\mu$-stability
  of $F$ yields the contradiction
  $\beta=\mu(F)<\mu(E')\leq\mu_{\rm max}(E)\leq\beta$.

  If $\varphi=0$, then
  $\kh^{-1}(E^\bullet)\cong\kh^0(G^\bullet)$ and ($\kh^{-1}(G^\bullet)\cong F$).
  Since the only common object of $\kt$ and $\kf$ is the trivial sheaf, the
  latter is only possible if $E\cong 0$. Hence, $E^\bullet\cong0$.

  If $\varphi$ is injective, then $\kh^{-1}(G^\bullet)=0$, i.e.\
  $G^\bullet\cong G:=\kh^0(G^\bullet)$ and we get a short exact
  sequence
  $$\xymatrix{0\ar[r]&F\ar[r]&E\ar[r]&G\ar[r]&0}$$
  in $\Coh(X)$ with $E\in\kf$ and $G\in\kt$.
  As $\mu_{\rm max}(E)\leq\beta$, the sheaf $G$ must be torsion.
  If $G$ is not concentrated in dimension zero, then $\deg(G)>0$
  and hence $\mu(E)=(\deg(G)+\deg(F))/\rk(F)>\mu(F)=\beta$
  contradicting $\mu_{\rm max}(E)\leq\beta$. If $G$ is
  concentrated in dimension zero, then
  $\Ext^1(G,F)\cong\Ext^1(F,G)^*\cong H^1(X,F\dual\otimes G)^*=0$,
  for $F$ is locally free. Thus, $E\cong F\oplus G$, which for
  $G\ne0$ contradicts the torsion freeness of $E$.

\smallskip

  {\bf iii) That's all.}
  Suppose $F^\bullet\in\ka$ is minimal and $F^\bullet\not\cong k(x)$
  for all closed points $x\in X$. The first two of the following
  claims, which we include here for completeness sake,
  correspond to b) of Lemma 6.1 in \cite{BK3}.

  {\it Claim 1.} $\kh^0(F^\bullet)=0$.\\
  Otherwise, there exists a surjection
  $\kh^0(F^\bullet)\twoheadrightarrow k(x)$ in $\Coh(X)$ and hence
  a non-trivial morphism $F^\bullet\to k(x)$ in
  $\ka$. As both objects are minimal, it would necessarily be an isomorphism.

  Hence, $F^\bullet\cong F[1]$ with $F\in\kf$.

  {\it Claim 2.} $F$ is locally free.\\
  If not, then there exists a short exact sequence
  $$\xymatrix{0\ar[r]&F\ar[r]&F'\ar[r]&k(x)\ar[r]&0}$$
  in $\Coh(X)$ with $F'$ still torsion free. Hence, $F'\in\kf$.
  Thus, the induced distinguished triangle
  $$\xymatrix{k(x)\ar[r]&F[1]\ar[r]&F'[1]}$$ yields a short exact
  sequence in $\ka$, contradicting the minimality of $F[1]$.

  {\it Claim 3.} $F$ is $\mu$-stable.\\
  If not, then there exists a short exact sequence
  \begin{equation}\label{ses3}
  \xymatrix{0\ar[r]&F_1\ar[r]&F\ar[r]&F_2\ar[r]&0}
  \end{equation}
  with $F_1,F_2$ torsion free, non-trivial and such that
  $\mu(F_1)\geq\mu(F)\geq\mu(F_2)$. Since $F\in\kf$, also
  $F_1,F_2\in\kf$. Therefore the shift of (\ref{ses3}) yields a
  short exact sequence in $\ka$ contradicting the minimality of
  $F^\bullet$.

  {\it Claim 4.} $\mu(F)=\beta$.\\
  Here we use Proposition \ref{ExSES} which shall be proved in the next
  section. It asserts that as soon as $\mu(F)<\beta$, there exists a
  short exact sequence
  $$\xymatrix{0\ar[r]&F\ar[r]&E\ar[r]&G\ar[r]&0}$$
  of $\mu$-stable vector bundles with $\mu(E)\leq \beta<\mu(G)$.
  The induced distinguished triangle
  $$\xymatrix{G\ar[r]&F[1]\ar[r]&E[1]}$$
  is then a short exact sequence in $\ka$ which again contradicts the
  minimality of $F[1]$. For an argument that does not make
  use of the full Proposition \ref{ExSES}, whose proof is
  unpleasantly long, see Remark \ref{BridgeRem}.
\end{proof}

In the `irrational' situation, the category
$\ka(\exp(B+i\omega))$ has the same minimal objects as $\Coh(X)$.
More precisely, the proposition yields:

\begin{corollary}
  Suppose the complex polarization $B+i\omega$ is chosen such that
  $\omega\in{\rm NS}(X)_\IQ$ and $(B.\omega)\not\in\IQ$. Then the
  only minimal objects in $\ka(\exp(B+i\omega))$ are the point
  sheaves $k(x)$.
\end{corollary}

\begin{proof}
If $\omega$ is rational, then for any sheaf $F$ the slope
$\mu(F)\in\IQ$. In particu\-lar, there are no $\mu$-stable vector
bundles with $\mu(F)=(B.\omega)$.
\end{proof}

\begin{remark} In the rational case, i.e.\ $B,\omega\in{\rm
NS}(X)_\IQ$, the minimal objects described by the proposition
share the property of the minimal objects in $\Coh(X)$ alluded to
in Example \ref{exCoh}: Every non-trivial object in
$\ka(\exp(B+i\omega))$ admits a surjection onto a minimal object.

Suppose $F^\bullet\in\ka$ with $\kh^0(F^\bullet)\ne0$. Then, any
surjection $\kh^0(F^\bullet)\twoheadrightarrow k(x)$ induces a
surjection $F^\bullet\twoheadrightarrow k(x)$ in $\ka$. If
$\kh^0(F^\bullet)=0$, then $F^\bullet\cong F[1]$. We may assume
that $F$ is $\mu$-stable and locally free. If
$\mu(F)=(B.\omega)$, then $F[1]\in\ka$ is minimal. If not, one
uses  Remark \ref{TechRem}, iii), which says that there always
exists a short exact sequence $0\to G\to F[1]\to E[1]\to0$ in
$\ka$ with $E$ $\mu$-stable, locally free and such that
$\mu(E)=\beta$.
\end{remark}

\section{Stable extensions: A technical fact}

 Let us fix a K\"ahler class $\omega$ on a projective K3 surface $X$ and
consider degree $\deg$ and slope $\mu$ with respect to $\omega$.

\begin{proposition}\label{ExSES}
  Fix $\beta\in \IR$. If $F$ is a $\mu$-stable vector bundle on
  $X$ with
  \begin{equation}\label{INEQ1}
  \mu(F)<\beta,
  \end{equation}
  then there exists a short exact sequence
  \begin{equation}\label{SES}
  \xymatrix{0\ar[r]&F\ar[r]&E\ar[r]&G\ar[r]&0}
  \end{equation}
  of $\mu$-stable vector bundles with
  \begin{equation}\label{INEQ2}
  \mu(E)\leq \beta < \mu(G).
  \end{equation}
\end{proposition}

The inequalities in (\ref{INEQ2}) impose numerical conditions on
the line bundles obtained as the determinants of the bundles in
(\ref{SES}). The existence of the line bundles is shown first.

In the following we let $L:=\det(F)$ and $r:=\rk(F)$.

\begin{lemma}\label{ExSESLB}
  Suppose $\mu(F)<\beta$. Then there exist a line bundle $L'$ and
  an integer $r'>0$ such that
  \begin{equation}\label{INEQ3}
  \frac{\deg(L)+\deg(L')}{r+r'}\leq\beta<\frac{\deg(L')}{r'}.
  \end{equation}
  Moreover, $L'$ can be chosen such that $L$ and $L'$ are linearly
  dependent.
\end{lemma}

\begin{proof}
  Suppose we can prove the existence of $L'$ and $r'$ satisfying
  (\ref{INEQ3}) for a twist $F\otimes H$ of $F$ by some (ample)
  line bundle $H$. Then $L'\otimes H^{-r'}$ would satisfy
  (\ref{INEQ3}) for $F$ itself. Note that $\mu(F)<\beta$ if and
  only if $\mu(F\otimes H)<\beta+\deg(H)$. Therefore, we may
  add the simplifying assumption $0<\deg(L)$.

  The line bundle $L'$ will be chosen of the form $L':= L^{\ell'}$
  for some integer $\ell'>0$. Dividing (\ref{INEQ1}) and
  (\ref{INEQ3}) by $\deg(L)$ the two inequalities become (with
  $\tilde\beta:=\beta/\deg(L)$):
  $$1<r\tilde\beta~~{\rm
  and}~~\frac{1+\ell'}{r+r'}\leq\tilde\beta<\frac{\ell'}{r'}.$$
  The latter is equivalent to
  $$\tilde\beta<\frac{\ell'}{r'}\leq\tilde\beta+\frac{\tilde\beta
  r-1}{r'}.$$

  To conclude recall the standard fact that for any $x\in\IR_{>0}$
  and all $\varepsilon >0$ there exists a rational number
  $\frac{a}{b}$ with $a,b$ positive integers such that
  $0<\frac{a}{b}-x<\frac{\varepsilon}{b}$. Apply this to
  $x:=\tilde \beta$ and $\varepsilon:=\tilde\beta r-1$ and set
  $\ell'=a$ and $r':=b$.
\end{proof}

\begin{remark}\label{BridgeRem} Although
of interest in its own right, Proposition \ref{ExSES} is in this
note only used in the proof of Proposition \ref{classmin}. As the
proof of Proposition \ref{ExSES} is rather lengthy, I'm very
grateful to  Tom Bridgeland who pointed out the following
shortcut to the argument in {\it Claim 4} in Section
\ref{SectMinObj}.

Suppose $L'$ and $r'$ are as in the Lemma and suppose there exists
a $\mu$-stable vector bundle $E$ with $\chi(F,E)>0$. Then
$\Hom(F,E)\ne0$ or $\Hom(E,F)\ne0$, but the latter is excluded by
stability. A non-trivial morphism $F\to E$ with $E$ $\mu$-stable
of slope $\leq \beta$ gives rise to a non-trivial morphism
$F[1]\to E[1]$ in $\ka$. But in the proof of Proposition
\ref{classmin} the object $F[1]$ was supposed to be simple and,
therefore, any non-trivial morphism from $F[1]$ is injective in
$\ka$. However, the quotient $G^\bullet$ (in $\ka$!) of $F[1]\to
E[1]$ satisfies $({\rm
c}_1(G^\bullet)-\rk(G^\bullet)B.\omega)=-(\ell'.\omega)+r'\beta<0$,
which contradicts $G^\bullet\in \ka$.

 Finally, one observes that
$\chi(F,E)=-\langle v(F),v(E)\rangle$, but if
$v(E)=(r+r',\ell+\ell',s')$, then $\langle
v(F),v(E)\rangle=(\ell.\ell+\ell')-rs'-(r+r')s<0$ for $s'\gg0$.
Thus, $\mu$-stable vector bundles $E$ of rank $r+r'$  with
$\det(E)\cong L\otimes L'$ and ${\rm c}_2(E)\gg0$, which
certainly exist, will yield the above contradiction to the
minimality of $F[1]$.
\end{remark}

Let us now prepare the proof of Proposition \ref{ExSES}. In the
course of the proof we shall make use of the existence of
$e$-stable vector bundles. Here $e$ is a real number, usually
positive, and a torsion free sheaf $G$ is called $e$-stable if
for all subsheaves $0\ne G_1\subset G$ with $\rk(G_1)<\rk(G)$ one
has
$$ \mu(G_1)<\mu(G)-\frac{e}{\rk(G_1)}.$$
For $e>0$ this is in general a stronger version of
$\mu$-stability.

O'Grady proved the existence of $e$-stable vector bundles with
large second Chern number: For fixed $e$, $L'$, and $r'$ and
$c\gg0$, there exists an $e$-stable vector bundle $G$ with
$\det(G)\cong L'$, $\rk(G)=r'$, and ${\rm c}_2(G)=c$. The bound
can be made effective (see \cite{OG} or \cite[Thm.\ 9.11]{HL}).

\bigskip

\bigskip

{\it Proof of Proposition \ref{ExSES}.}
  In the following we let $L'$ and $r'$ be as in Lemma
  \ref{ExSESLB} and we assume furthermore $r'\geq r$. The vector
  bundle $G$ will be chosen such that $\det(G)\cong L'$ and
  $\rk(G)=r'$. The remaining numerical invariant of $G$ is its
  second Chern number ${\rm c}_2(G)$, which will have to be chosen
  large enough.

  Consider any extension (possibly trivial)
  \begin{equation}\label{SES2}
  \xymatrix{0\ar[r]&F\ar[r]&E\ar[r]&G\ar[r]&0}
  \end{equation}
  and a proper saturated subbundle $0\ne E_1\subset E$. Then let
  $F_1:=F\cap E_1$ and $G_1:=E_1/F_1={\rm Im}(E_1\to G)$.

  The following short-hands will be used throughout:
  $\ell:=\deg(L)$, $\ell':=\deg(L')$, $\ell_1:=\deg(F_1)$,
  $\ell_1':=\deg(G_1)$, $r_1:=\rk(F_1)$, and $r_1':=\rk(G_1)$.

\smallskip

  {\bf i)} If $G_1=0$, then $F_1=E_1$. Hence, by the $\mu$-stability of
  $F$ one finds
  $$\mu(E_1)=\mu(F_1)\leq\mu(F)<\mu(E).$$

\smallskip

  {\bf ii)} Suppose $0<\rk(G_1)<\rk(G)$.\\
  {\it Claim.} If $G$ is $e$-stable for some
  $e\geq\left(\frac{r'r-r}{r+r'}\right)\left(\frac{\ell'}{r'}-\frac{\ell}{r}\right)$,
  then $\mu(E_1)<\mu(E)$.

  The $\mu$-stability of $F$ yields
  \begin{eqnarray*}
  \mu(E_1)&=&\frac{\ell_1+\ell_1'}{r_1+r_1'}=\frac{\ell_1}{r_1}\cdot\frac{r_1}{(r_1+r_1')}
  +\frac{\ell'_1}{r'_1}\cdot\frac{r'_1}{(r_1+r_1')}\\
  &\leq&\frac{\ell}{r}\cdot\frac{r_1}{(r_1+r_1')}+\frac{\ell_1'}{r_1'}\cdot\frac{r_1'}{(r_1+r_1')}
  \end{eqnarray*} (We leave it to the reader to verify that this makes
  sense also in the case $r_1=0$.)  Thus, it suffices to show that the right hand side is smaller
  than $\mu(E)=\frac{\ell+\ell'}{r+r'}$ or, equivalently, that
  $$\frac{\ell_1'}{r_1'}<\frac{(\ell+\ell')}{r_1'}\cdot
  \frac{(r_1+r_1')}{(r+r')}-\frac{\ell}{r}\cdot\frac{r_1}{r_1'}.$$
  Since $G$ is $e$-stable and $\rk(G_1)<\rk(G)$, one has
  $$\frac{\ell_1'}{r_1'}<\frac{\ell'}{r'}-\frac{e}{r_1'}.$$
  Hence $\mu(E_1)<\mu(E)$ if
  $$\frac{\ell'}{r'}-\frac{e}{r_1'}<\frac{(\ell+\ell')}{r_1'}
  \cdot\frac{(r_1+r_1')}{(r+r')}-\frac{\ell}{r}\cdot\frac{r_1}{r_1'}$$
  or, equivalently, if
  $$\left(\frac{\ell'}{r'}-\frac{\ell}{r}\right)
  \left(\frac{r_1'r-r'r_1}{r+r'}\right)<e.$$
  The maximum of the left hand side with $0\leq r_1\leq r$ and
  $0<r_1'<r'$ is attained for $r_1=0$ and $r_1'=r'-1$  and yields
  $\left(\frac{r'r-r}{r+r'}\right)\left(\frac{\ell'}{r'}-\frac{\ell}{r}\right)$.

\smallskip

  {\bf iii)} Consider now the remaining case that $\rk(G_1)=\rk(G)$. Since
  $E_1\subset E$ is a saturated and proper subsheaf, $F/F_1$ is torsion free and
  $\rk(F_1)<\rk(F)$. If $E_1$ is $\mu$-destabilizing, then
  $\mu(E_1)\geq\mu(E)$ or, equivalently,
  \begin{equation}\label{ineqpf}
  \ell_1\geq(\ell+\ell')\frac{r_1+r'}{r+r'}-\ell_1'.
  \end{equation}
  Consider $${\mathcal S}:=\{F_1\subset F~|~F/F_1~{\rm
  torsion~free},~\mu(F_1)\geq C\}$$
  with $C:=(\ell+\ell')\frac{r_1+r'}{r+r'}-\ell'$. Note that $C$
  is bigger than the right hand side of (\ref{ineqpf}).

  If $\xi\in \Ext^1(G,F)$ denotes the extension class of
  (\ref{SES2}) and $\eta\in \Ext^1(S,F/F_1)$ the extension class
  of the induced short exact sequence
  $$\xymatrix{0\ar[r]&F/F_1\ar[r]&E/E_1\ar[r]&S:=G/G_1\ar[r]&0,}$$
  then they yield identical classes in $\Ext^1(G,F/F_1)$ under the
  natural maps $\Ext^1(G,F)\to \Ext^1(G,F/F_1)$ respectively
  $\Ext^1(S,F/F_1)\to\Ext^1(G,F/F_1)$.

  Suppose we can choose $G$ and $\xi\in \Ext^1(G,F)$ such that:

  \begin{itemize}\item[(A)]
  For all $F_1\in{\mathcal S}$ the image of $\xi$ is not
  contained in the image of $\Ext^1(S,F/F_1)\to\Ext^1(G,F/F_1)$
  for any torsion quotient $G\twoheadrightarrow S$.\end{itemize}

  Then $E$ does not admit a destabilizing subsheaf with
  $\rk(G_1)=\rk(G)$.

\smallskip

  {\bf iv)} To conclude the proof it suffices to show that there exists
  an $e$-stable vector bundle $G$ with $\det(G)=L'$, $\rk(G)=r'$,
  $e$ as in ii), and an extension class $\xi\in\Ext^1(G,F)$
  satisfying (A).

  First note that due to a lemma of Grothendieck (see \cite[Lemma
  1.7.9]{HL}) the family ${\mathcal S}$ is bounded. Next consider the
  reflexive hull of $F/F_1$ which sits in a short exact sequence
  $$\xymatrix{0\ar[r]&F/F_1\ar[r]&(F/F_1)\ddual\ar[r]&T\ar[r]&0}$$
  for some torsion sheaf $T$. Applying $\Hom(S,~~)$ yields a
  bijection
  $\Hom(S,T)\cong\Ext^1(S,F/F_1)$ for any torsion quotient
  $G\twoheadrightarrow S$. Thus there is a commutative diagram
  $$\xymatrix{\Hom(S,T)\ar@{^{(}->}[d]\ar[r]^-\sim&\Ext^1(S,F/F_1)\ar[d]\\
  \Hom(G,T)\ar[r]&\Ext^1(G,F/F_1),}$$ which shows that the image of
  $\Ext^1(S,F/F_1)\to\Ext^1(G,F/F_1)$ is contained in the image of
  $\Hom(G,T)\to\Ext^1(G,F/F_1)$ which is independent of $S$.

  Since ${\mathcal S}$ is a bounded family, the length of the occurring
  sheaves $T=(F/F_1)\ddual/(F/F_1)$ with $F_1\in{\mathcal S}$ is
  bounded by say $a$. Hence, for any $F_1\in{\mathcal S}$ the
  dimension of \begin{eqnarray*}
  &\bigcup_{G\twoheadrightarrow S}{\rm
  Im}\left(\Ext^1(S,F/F_1)\to\Ext^1(G,F/F_1)\right)\\&
  \subset V_{F_1}:={\rm Im}\left(\Hom(G,T)\to\Ext^1(G,F/F_1)\right)
  \end{eqnarray*}
  is bounded by $r'a$.

  Consider the pre-image $V'_{F_1}\subset\Ext^1(G,F)$ of
  $V_{F_1}\subset \Ext^1(G,F/F_1)$ under
  $$\Ext^1(G,F)\to\Ext^1(G,F/F_1).$$
  In order to find $\xi$ satisfying (A) it suffices
  to show that for ${\rm c}_2(G)\gg0$ and $G$ generic, the
  algebraic set $\bigcup_{F_1\in{\mathcal S}} V'_{F_1}\subset
  \Ext^1(G,F)$ has dimension strictly smaller than $\ext^1(G,F)$.

  By definition, $\dim (V_{F_1}')\leq\dim (V_{F_1})+\ext^1(G,F_1)$.
  Thus, $$\dim(\bigcup_{F_1\in{\mathcal S}}V_{F_1}')\leq
  r'a+\sup\{\ext^1(G,F_1)~|~F_1\in{\mathcal S}\}+\dim({\mathcal S}).$$

  Write
  $\ext^1(G,F_1)=\hom(G,F_1)+\hom(F_1,G)-\chi(G,F_1)$. Clearly, if
  $G$ is $\mu$-stable, then $\Hom(G,F_1)\subset\Hom(G,F)=0$.
  Since ${\mathcal S}$ is bounded, there exists a constant $\mu_0$
  such that $\mu_0\leq\mu(\hat F_1)$ for all quotients $F_1\twoheadrightarrow
  \hat F_1$ of a sheaf  $F_1\in{\mathcal S}$.

  If $G$ is $e$-stable with
  $e\geq\mu(G)-\mu_0=\frac{\ell'}{r'}-\mu_0$, then
  $\Hom(F_1,G)=0$. Indeed, for the image $\hat F_1$ of a
  non-trivial morphism $F_1\to G$ one would obtain the
  contradiction $\mu_0\leq\mu(\hat F_1)<\mu(G)-e\leq\mu_0$.
  (Recall $r'\geq r$.)

  The Riemann--Roch formula then shows that for $e\geq\frac{\ell'}{r'}-\mu_0$ and an
  $e$-stable vector bundle $G$ the dimension
  $\ext^1(G,F_1)=-\chi(G,F_1)$ grows like $r_1{\rm c}_2(G)$ for
  ${\rm c}_2(G)\to\infty$. On the other hand, $\ext^1(G,F)$ grows
  at least like $r{\rm c}_2(G)$. Hence,
  $$\dim(\bigcup_{F_1\in{\mathcal S}}V'_{F_1})<\dim\Ext^1(G,F)$$ for
  ${\rm c}_2(G)\gg0$ and $G$ an $e$-stable vector bundle
  for any $e\geq\frac{\ell'}{r'}-\mu_0$.

  Eventually choose
  $e\geq\max\{\left(\frac{rr'-r}{r+r'}\right)
  \left(\frac{\ell'}{r'}-\frac{\ell}{r}\right),
  \frac{\ell'}{r'}-\mu_0\}$.
  \qed
\begin{remark}\label{TechRem}
  i) The proof also shows that $G$ (and hence $E$) can be chosen such
  that $\det(F)$ and $\det(G)$ are linearly dependent.

  ii) The proposition holds true for any torsion free $\mu$-stable
  sheaf $F$. Again, $G$ can be chosen locally free, but $E$ would
  only be torsion free in this more general situation.

  iii) One can be a bit more specific about the slope $\mu(E)$. In
  fact, any slope that could in principle be realized can
  also be realized as a slope $\mu(E)$. For example if $\omega$ and
  $\beta$ are both rational, then we can find $E$ such that
  $\mu(E)=\beta$.
\end{remark}

\section{FM-partners via $\mu$-stable vector bundles}

Due to results of Mukai \cite{Mu} and Orlov \cite{Or}, one knows
that two K3 surfaces $X$ and $X'$ are derived equivalent if and
only if $X'$ is a fine moduli space of stable sheaves on $X$ (see
also \cite[Ch.\ 10]{HFM}). A priori `stable' in this context
means `Gieseker stable' (and not $\mu$-stable) and the sheaves are
just torsion free (and not locally free). However, as will be
shown in this section, the stronger result holds true, i.e.\ one
can work with $\mu$-stable locally free sheaves. The result might
be known to the experts -- the techniques certainly are. In
particular, Yoshioka treats this question in \cite[Lemma 2.1]{Y2},
but I was not always absolutely sure about the assumptions in
\cite{Y2} and in the article \cite{Y3} it is based on. In any
case, as the explicit statement, crucial for the rest of the
paper, does not seem to be in the literature and for the reader's
convenience, we include a complete proof here.

\begin{proposition}\label{FMviamu}
Two K3 surfaces $X$ and $X'$ are derived equivalent if and only if
either $X\cong X'$ or $X'$ is isomorphic to a fine moduli space
of $\mu$-stable vector bundles on $X$.
\end{proposition}

\begin{proof}
One direction is a special case of Mukai's result. So, we only
have to prove that if $\Db(X')\cong \Db(X)$, then either $X\cong
X'$ or $X'$ is isomorphic to a fine moduli space of $\mu$-stable
vector bundles on $X$.

i) We shall often need the following fact: If $\Phi:\Db(Y)\congpf
\Db(Y')$ is an equivalence between two K3 surfaces with
$\Phi^H(0,0,1)=(0,0,1)$, then $Y\cong Y'$. Indeed, $\Phi^H$ then
induces a Hodge isometry
$(0,0,1)_Y^\perp\congpf(0,0,1)_{Y'}^\perp$ and a Hodge isometry
of the quotients
$$H^2(Y,\IZ)\cong(0,0,1)^\perp_Y/(0,0,1)\IZ\congpf(0,0,1)^\perp_{Y'}/(0,0,1)\IZ\cong
H^2(Y',\IZ).$$ By the Global Torelli theorem this implies $Y\cong
Y'$.

ii) Orlov's proof \cite{Or} (or \cite[Sect.\ 10.2]{HFM}) shows
that if $\Db(X)\cong\Db(X')$, then $X'$ is isomorphic to a moduli
space $M_H(v)$ of Gieseker stable (with respect to $H$) sheaves on
$X$ with Mukai vector $v=(r,a\ell,s)$, where
\begin{equation}\label{crucform}
\ell\in{\rm NS}(X)~{\rm primitive},~{\rm
g.c.d}(r,a(\ell.H),s)=1,~{\rm and} ~a^2(\ell.\ell)=2rs.
\end{equation}
 Clearly, the last equality is
just expressing $\dim (M_H(v))=\dim(X')=2$. The fact that $r,a,s$
are coprime not only ensures that every Gieseker semistable sheaf
is Gieseker stable, but also that the moduli space is fine (see
\cite[Cor.\ 4.6.7]{HL}). Moreover, one may assume $r\geq2$ if $X$
and $X'$ are not already isomorphic.

iii) Observe that if ${\rm g.c.d.}(r,a)=1$ and if $H'$ is a
polarization not lying on a wall, then
$M_{H'}(v)=M_{H'}(v)^{\mu{\rm ss}}$, i.e.\ every Gieseker
$H'$-stable sheaf with Mukai vector $v$ is $\mu$-stable (see
\cite[Thm.\ 4.C.3]{HL}). Furthermore, if $\ke$ and $\ke'$ are the
universal families on $X'\times X=M_H(v)\times X$ and
$M_{H'}(v)\times X$, respectively, then the equivalence
$$\xymatrix@=27pt{\Db(X')\ar[r]^-\sim_-{\Phi_\ke}&\Db(X)\ar[r]^-\sim_-{\Phi_{\ke'}^{-1}}&\Db(M_{H'}(v))}$$
sends $(0,0,1)$ to ${\Phi_{\ke'}^H}^{-1}(\Phi_\ke^H(0,0,1))
={\Phi_{\ke'}^H}^{-1}(v)=(0,0,1)$. Hence, due to i), one has
$X'\cong M_{H'}(v)$, i.e.\ $X'$ is isomorphic to a fine moduli
space of $\mu$-stable sheaves.

iv) Let us now treat the case of Picard number one. The argument
used here is not very geometric, as it makes use of a counting of
primitive embeddings into the K3 lattice. (Is there a better one?)

According to \cite[Thm.\ 2.1]{HLOY} every K3 surface $X'$ derived
equivalent to a given K3 surface $X$ with $\Pic(X)=\IZ\ell$ is
isomorphic to a moduli space $M_H(v)$ with $v=(r,\ell,s)$ and
${\rm g.c.d.}(r,s)=1$. In particular, the determinant is primitive
and this ensures that $\mu$-semistability implies $\mu$-stability
(use iii) or, more directly, the assumption $\rho(X)=1$). The
proof of this result relies in an explicit counting of all
Fourier--Mukai partners in \cite{Ogui} and a counting of all
Fourier--Mukai partners arising as one of these moduli spaces
(see also \cite{St}).

v) For the case of $\rho(X)\geq2$ we need to modify a given Mukai
vector by spherical twists and line bundle twists. This shall be
prepared now.

Recall that $\ko\in\Db(X)$ is a spherical object and therefore
induces an autoequivalence $T_\ko:\Db(X)\congpf\Db(X)$, the
spherical twist. Its action $T^H_\ko$ on $\widetilde H(X,\IZ)$
interchanges the generators of $H^0$ and $H^4$ (up to sign) and
leaves invariant $H^2(X,\IZ)$. If $v$ satisfies (\ref{crucform}),
then $v':=\pm T^H_\ko(v)$ does as well. (Choose the sign such that
the rank of $v'$ is non-negative.) Thus, the moduli space
$M_H(v')$ is non-empty and fine. If $\ke'$ denotes the universal
family on $M_H(v)\times X$, then the cohomological Fourier--Mukai
transform induced by the composition
$$\xymatrix@=27pt{\Db(M_H(v'))\ar[r]^-\sim_-{\Phi_{\ke'}}&\Db(X)\ar[r]^-\sim_{T^{-1}_\ko}&\Db(X)
\ar[r]^-\sim_-{\Phi^{-1}_\ke}&\Db(M_H(v))}$$ maps $(0,0,1)$ to
$(0,0,1)$. Hence, due to i), one has  $M_H(v')\cong M_H(v)$. (If
$v'=-T_\ko^H(v)$, then replace $T_\ko^{-1}$ by $T_\ko^{-1}[1]$.)

If $\tilde L\in\Pic(X)$, then the autoequivalence $\tilde
L\otimes(~~):\Db(X)\congpf\Db(X)$ maps (on the cohomology)
$v\in\widetilde H(X,\IZ)$ to $\tilde v:=\exp({\rm c}_1(\tilde
L))\cdot v$, which once again satisfies (\ref{crucform}). As
above,
$$\xymatrix@=27pt{\Db(M_H(\tilde v))\ar[r]^-\sim_-{\Phi_{\tilde\ke}}
&\Db(X)\ar[r]^-\sim_{\tilde L^*\otimes(~~)}&\Db(X)
\ar[r]^-\sim_-{\Phi^{-1}_\ke}&\Db(M_H(v))}$$ sends $(0,0,1)$ to
$(0,0,1)$. Hence, $M_H(\tilde v)\cong M_H(v)$.

Summarizing, we conclude that the Mukai vector $v$ can be modified
at will by $T_\ko^H$ and $\exp(\tilde \ell)$ with
$\tilde\ell\in{\rm NS}(X)$ without changing the isomorphism type
of the moduli space.

vi) We are now ready to treat the case $\rho(X)\geq2$. Suppose
$v=(r,a\ell,s)$ satisfies (\ref{crucform}) and
$r+a\ell=\alpha(r'+a'\ell)$ with ${\rm g.c.d.}(r,a)=1$ (or,
equivalently, $r'+a'\ell\in (H^0\oplus H^2)(X,\IZ)$ is primitive).
Since $\rho(X)\geq2$, we may choose a primitive
$\tilde\ell\in{\rm NS}(X)$ linearly independent of $\ell$.
Consider $$\exp(\tilde\ell)\cdot
v=\left(r,a\ell+r\tilde\ell,\tilde
s:=s+r(\tilde\ell.\tilde\ell)/2+a(\tilde\ell.\ell)\right).$$ Then
$r+a\ell+r\tilde\ell=\alpha(r'+a'\ell+r'\tilde\ell)$ and
$r'+a'\ell+r'\tilde\ell$ primitive. Furthermore, as $\alpha$
divides $r$ and $a$ and ${\rm g.c.d.}(r,a,s)=1$, one has ${\rm
g.c.d.}(\alpha,\tilde s)=1$.

Hence, for $T_\ko^H(\exp(\tilde\ell)\cdot v)$ rank and determinant
are coprime. Thus, either $X\cong X'$ or $X'$ is isomorphic to a
moduli space of $\mu$-stable sheaves on $X$.

vii) In the last step, one has to show that either $X'\cong X$ or
$X'$ is isomorphic to a fine moduli space of $\mu$-stable locally
free(!) sheaves on $X$. According to iv) and vi), $X\cong X'$ or
$X'$ is isomorphic to a fine moduli space $M_H(v)$ of $\mu$-stable
sheaves with $v=(r,\ell,s)$ and $r>0$. Suppose there exists a
$\mu$-stable torsion free sheaf $[E]\in M_H(v)$ which is not
locally free.

Consider the natural sequence
$$\xymatrix{0\ar[r]&E\ar[r]&E\ddual\ar[r]&S\ar[r]&0,}$$
for which $S$ is torsion. Since $E$ is $\mu$-stable  and $S$ is
concentrated in dimension zero, $E\ddual$ is a $\mu$-stable
vector bundle. As Mukai observed in \cite[Prop.\ 3.9]{Mu}, the
assumption that $[E]\in M_H(v)$ and hence $\ext^1(E,E)=2$ implies
$S\cong k(x)$. Deforming $x\in X$ and the surjection
$E\ddual\twoheadrightarrow k(x)$ with it, yields a family of
kernels of dimension $\dim(X)+\dim(\IP(E\dual(x)))=2+r-1$. Hence,
$r=1$, which implies $X\cong X'$.

Thus all sheaves $[E]\in M_H(v)$ are locally free.
\end{proof}

\section{Derived equivalence yields abelian
equivalence}\label{SectDA}

Suppose $X$ and $X'$ are two derived equivalent K3 surfaces. By
Proposition \ref{FMviamu},  there exists an isomorphism $X\cong
X'$, the trivial case that will not be discussed, or $X'$ is
isomorphic to a fine moduli space of $\mu$-stable vector bundles
on $X$. In the latter case we shall denote the universal family on
$X\times X'$ by $\ke$. So, for any closed point $y\in X'$ the
restriction $\ke_y:=\ke|_{X\times\{y\}}$ is a $\mu$-stable vector
bundle on $X$. Here, $\mu$-stability is taken with respect to a
K\"ahler class $\omega\in{\rm NS}(X)_\IR$, which can be chosen
integral and such that $(\omega.\omega)>2$. Clearly, the slope
$\mu(\ke_y)$ is independent of $y\in Y$ and there exists a
rational B-field $B\in{\rm NS}(X)_\IQ$ such that
$\mu(\ke_y)=(B.\omega)=:\beta$ for all $y\in X'$.

In the following we shall consider the $t$-structure induced by
$\omega$ and $B$ (respectively $\beta$) and consider its heart
$$\ka:=\ka(\exp(B+i\omega))=\ka(\beta).$$

\begin{remark}
  Bridgeland shows that
  $Z(F^\bullet)=\langle\exp(B+i\omega),v(F^\bullet)\rangle$ is a
  stability function that satisfies the HN-property (see
  \cite[Sect.\ 6,9]{BK3}). So, the $t$-structure induced by
  $\exp(B+i\omega)$ and together with $Z$ define a stability condition
  on $\Db(X)$.

  We shall only use that
  $Z(F^\bullet)\in\IR_{>0}\exp(i\pi\phi(F^\bullet))$ with
  $0<\phi(F^\bullet)\leq1$ (see Remark \ref{RemStabBrid})
  for any $0\ne F^\bullet\in\ka$.
\end{remark}

As was proved already by Mukai in \cite{Mu}, the Fourier--Mukai
transform
$$\xymatrix{\Phi:=\Phi_{\ke[1]}:\Db(X')\ar[r]^-\sim&\Db(X)}$$
with kernel $\ke[1]\in\Db(X'\times X)$ is an exact equivalence.
Note that for a closed point $y\in X'$ one has
$$\Phi(k(y))\cong \ke_y[1].$$
Furthermore, its inverse $\Psi:=\Phi^{-1}$ is the Fourier--Mukai
transform $\Phi_{\ke\dual[1]}$ with kernel
$\ke\dual[1]\in\Db(X\times X')$.

The stability condition given by $(\ka,Z)$ induces via the
equivalence
$$\xymatrix{\Psi:\Db(X)\ar[r]^-\sim&\Db(X')}$$
a stability condition on $\Db(X')$. More precisely, one sets
$$\ka':=\Psi(\ka)$$ and $$Z'(F^\bullet):=\left\langle
\Psi^H(\exp(B+i\omega)),v(F^\bullet)\right\rangle=
\left\langle\exp(B+i\omega),v(\Phi^H(F^\bullet))\right\rangle.$$
Here, $\Psi^H:\widetilde H(X,\IZ)\congpf \widetilde H(X',\IZ)$ and
its inverse $\Phi^H:\widetilde H(X',\IZ)\congpf \widetilde
H(X,\IZ)$ are the naturally induced Hodge isometries.

Clearly, ${\rm Im}(Z'(F^\bullet))\geq0$ for any
$F^\bullet\in\ka'$ and  any $k(y)=\Psi(\ke_y[1])$ satisfies
$Z'(k(y))=Z(\ke_y[1])\in\IR$. As by Proposition \ref{classmin} the
$\mu$-stable vector bundles $\ke_y$ yield minimal objects
$\ke_y[1]\in\ka$, this shows that all point sheaves $k(y)$ are
minimal objects in $\ka'$, hence stable, of phase $\phi'(k(y))=1$.

Next recall that
$$\Psi^H(\exp(B+i\omega))=\lambda\exp(B'+i\omega')$$ for some ample
class $\omega'\in {\rm NS}(X')$ and a positive integer $\lambda$
(see \cite[Sect.\ 5]{HS} or \cite[Lemma 7.1]{Yos1}). (Note that
since $\Psi^H$ is an integral Hodge isometry, $\lambda$ and
$\omega'$ are indeed integral.)

We wish to compare the abelian category $\ka'\subset \Db(X')$
with the heart $\ka(\exp(B'+i\omega'))$ associated to the ample
class $\omega'$ and the B-field $B'$.

\begin{proposition}\label{indabequ}
  The two abelian subcategories $\ka'$ and $\ka(\exp(B'+i\omega'))$
  coincide, i.e.\
  $$\ka'=\ka(\exp(B'+i\omega'))\subset\Db(X').$$
\end{proposition}

\begin{proof}
  The proof is an immediate consequence of the discussion in
  Section 6 of \cite{BK3}. Bridgeland
  shows that the heart of any stability condition on $\Db(X')$ for
  which all point sheaves $k(y)$ are stable of phase one is of the
  type claimed by the proposition. For the convenience of the reader
  (and because the result in \cite{BK3} is not quite phrased as explicitly as
  our assertion) we include the argument. (Also, we don't really
  use that $\ka'$ is the heart of  a stability condition.)

  To shorten the notation, we let $\beta':=(B'.\omega')$ and denote
  the abelian category $\ka(\exp(B'+i\omega'))$ simply by
  $\ka(\beta')$ (the polarization $\omega'$ is understood).
  Similarly, the torsion pair of $\Coh(X')$ defining $\ka(\beta')$
  is denoted $(\kt(\beta'),\kf(\beta'))$.

  Bridgeland defines a torsion pair $(\kt',\kf')$  in $\Coh(X')$ by
  $$\kt':=\ka'\cap\Coh(X')~~{\rm and}~~\kf':=\ka'[-1]\cap\Coh(X'),$$ for
  which it is easy to verify that its tilt yields $\ka'$.
  All what is needed to prove this is collected in \cite[Lemma
  6.1]{BK3} (see also the arguments in the proof of
  Proposition \ref{classmin}):\\
  a) Any object $F^\bullet\in\ka'$ is concentrated in
  degree $0$ and $-1$ and  $\kh^{-1}(F^\bullet)$ is torsion free.\\
  b) If $F^\bullet\in\ka'$ is stable of phase one, then either
  $F^\bullet\cong k(y)$ or $F^\bullet\cong F[1]$ with $F$ locally
  free. (See the arguments in the proof of Proposition \ref{classmin}.)\\
  c) $\Coh(X')\subset\ka'\cup\ka'[-1]$.

  \smallskip

  Thus, it
  suffices to prove that $\kt'=\kt(\beta')$ and $\kf'=\kf(\beta')$
  or, equivalently, $\kt(\beta')\subset\kt$ and
  $\kf(\beta')\subset\kf'$.

  One first proves $\kt(\beta')\subset\kt'$. Let $F\in\kt(\beta')$.
  Due to a) and c) any torsion sheaf is contained in $\kt'$. Since $F/T(F)$
  is again in $\kt(\beta')$ and $\kt'$ is closed under extension,
  we can assume that $F$ is torsion free. Next consider the
  HN-filtration of $F$ with quotients $F_{i+1}/F_i$ which are all,
  due to the definition of $\kt(\beta')$,  torsion free and
  $\mu$-semistable with $\mu(F_{i+1}/F_i)>\beta'$. Again using that
  $\kt'$ is closed under extension, we can thus restrict to the case
  that $F$ is $\mu$-semistable and further, by using the Jordan--H\"older
  filtration, to $F$ $\mu$-stable.
  Now consider the decomposition of $F$ with respect to the torsion
  pair $(\kt',\kf')$, i.e.\ the short exact sequence
  $$\xymatrix{0\ar[r]&F_1\ar[r]&F\ar[r]&F_2\ar[r]&0}$$
  with $F_1\in \kt'$ and $F_2\in\kf'$. If both $F_1$ and $F_2$ are
  non-trivial, one obtains the contradiction
  $$\beta'\leq\mu(F_1)<\mu(F)<\mu(F_2)\leq\beta'.$$ Hence, either
  $F\in\kt'$ or $F\in\kf'$. The latter can be excluded as follows.
  Any object in $\kf'$ is of the form $\Psi(E^\bullet)[-1]$
  with $E^\bullet\in\ka$. Therefore, ${\rm Im}(Z'(F))=-{\rm Im}(Z(E^\bullet))\leq0$
  contradicting $F\in\kt(\beta')$, i.e.\  $\mu(F)>\beta'$.

  Next, one shows $\kf(\beta')\subset\kf'$. For this purpose pick
  $F\in\kf(\beta')$, which by definition is torsion free and such
  that all HN-factors have slope $\leq\beta'$. As above, it
  suffices to show that any $\mu$-stable torsion free sheaf $F$
  with $\mu(F)\leq\beta'$ is contained in $\kf'$. The same reasoning
  as above allows one to conclude that either $F\in\kf'$ or
  $F\in\kt'$. If $F\in\kt'$, then choose a point $y\in X'$ and a
  surjection $\varphi:F\twoheadrightarrow k(y)$ in $\Coh(X')$.
  Considered as a morphism in $\ka'$ it is non-trivial and, since
  $k(y)$ is minimal in $\ka'$, also surjective. Thus, its kernel $F_1$ is
  again contained in $\kt'\subset\ka'$ and satisfies
  $Z'(F_1)\in\IR_{<0}$. Hence $Z'(F)=Z'(k(y))+Z'(F_1)$.
  Now continue with $F_1$. This leads to $Z'(F)=kZ'(k(y))+
  \sum_{i=1}^kZ'(F_i)$ with $Z'(F_i)\in\IR_{<0}$ and hence to a
  contradiction for $k\to\infty$.
\end{proof}

By definition $\Psi$ induces an equivalence $\ka\congpf\ka'$. This
yields the second part of Theorem \ref{ThmfirstIntro}.

\begin{corollary}\label{CorDEyieldsAE}
If $X$ and $X'$ are two derived equivalent K3 surfaces, then
there exist complexified polarizations $B+i\omega$ and
$B'+i\omega'$ on $X$ respectively $X'$ and an equivalence
$$\ka_X(\exp(B+i\omega))\cong\ka_{X'}(\exp(B'+i\omega')).$$
\end{corollary}

\begin{remark}
Suppose $\Phi:\Db(X)\congpf\Db(X')$ is any equivalence such that
$\Phi^H(\exp(B+i\omega))=\exp(B'+i\omega')$ up to a positive
scalar, e.g.\ $\Phi$ induced by a universal sheaf of Gieseker
stable torsion free sheaves. Then $\Phi(\ka_X(\exp(B+i\omega)),
Z)$ is the heart of a stability condition on $X'$ in the same
fibre of $\pi:{\rm Stab}(X')\to\kp^+_0(X')$ as
$(\ka_{X'}(\exp(B+i\omega), Z')$. If it is contained in
Bridgeland's distinguished component $\Sigma(X')$, then there
exists an autoequivalence $\Psi$ of $\Db(X')$ such that
$\Psi(\Phi(\ka_X(\exp(B+i\omega))))=\ka_{X'}(\exp(B'+i\omega'))$.
The above calculation shows that for the Fourier--Mukai
equivalence induced by the universal family of $\mu$-stable
vector bundles we do not need to check whether the image is
contained in $\Sigma(X')$ and can set directly $\Psi={\rm id}$.
\end{remark}
\section{Twistor space interpretation}\label{SectTw}

Suppose $\omega$ is a K\"ahler class on a not necessarily
projective K3 surface $X$. Identifying $\omega$ with its unique
Ricci-flat (hyperk\"ahler) representative allows one to write
down an explicit model of the associated twistor space
$\pi:\kx\to\IP^1$. If $I$ is the complex structure defining $X$
and $J,K$ are the complementary ones determined by $\omega$, then
$\kx$ is $X\times\IP^1$ endowed with the complex structure acting
by $(\lambda_1I+\lambda_2J+\lambda_3K,I_{\IP^1})$ in
$(x,\lambda:=(\lambda_1+\lambda_2+\lambda_3))\in X\times\IP^1$.

In particular, $\pi$, which is by definition the second
projection, is holomorphic and for any $x\in X$ the curve
$L_x:=\{x\}\times\IP^1\subset X\times\IP^1=\kx$ is a holomorphic
section of $\pi$, the twistor sections.

\begin{itemize}
\item[i)] {\it To $L_x$ one associates $\ko_{L_x}\in\Coh(\kx)$, which
has the property that $\ko_{L_x}|_{X=\pi^{-1}(I)}\cong k(x)$.}
\end{itemize}

Suppose $E$ is a holomorphic vector bundle on $X$ which is
$\mu$-stable with respect to $\omega$. Then $E$ admits a unique
Hermite--Einstein metric. The curvature of the induced Chern
connection satisfies the Hermite--Einstein equation
$i\Lambda_\omega F_\nabla=\eta\cdot {\rm id}$ with
$\eta\int_X\omega^2=4\pi\mu(E)$. Thus, if $\mu(E)=0$, the
equation reads $i\Lambda_\omega F_\nabla=0$. The bundle $E$ can be
viewed as a complex vector bundle simultaneously on all the
fibres $\kx_\lambda=\pi^{-1}(\lambda)$ and, as was first observed
by Itoh (see e.g.\ \cite{B,V}), the $(0,1)$-part with respect to
$\lambda$ defines again a $\bar\partial$-operator. Thus,
$(E,\nabla^{(0,1)_\lambda})$ is a holomorphic vector bundle on
$\kx_\lambda$. Moreover, these bundles glue to a holomorphic
vector bundle on $\kx$.

\begin{itemize}
\item[ii)]{\it To any $\mu$-stable vector bundle $E$ on $X$ with
$\mu(E)=0$ one associates a distinguished vector bundle
$\ke\in\Coh(\kx)$ with $\ke|_X\cong E$.}
\end{itemize}

In the above situation consider a generic fibre $\kx_\lambda$ and
$F\in\Coh(\kx_\lambda)$. Since $\kx_\lambda$ does not contain any
curves, the torsion of $F$ is concentrated in dimension zero and
therefore admits a filtration with quotients isomorphic to some
$k(x)$.

The fibre $\kx_\lambda$ is endowed with a natural K\"ahler form
$\omega_\lambda=\lambda_1\omega_I+\lambda_2\omega_J+\lambda_3\omega_K$
and $\deg$ and $\mu$ are considered with respect to it. Due to
the generic choice of $\lambda$, one has $\deg(L)=0$ for any line
bundle $L\in\Pic(\kx_\lambda)$. In particular, the reflexive hull
$F\ddual$ is (if not trivial) a $\mu$-polystable vector bundle,
which can by the above procedure obtained as a restriction of the
direct sum of some $\ke$.

Thus, for generic $\lambda\in\IP^1$ the coherent sheaves
$$\ko_{L_x}|_{\kx_\lambda}~~{\rm and}~~\ke|_{\kx_\lambda}$$
with $\ko_{L_x}$ and $\ke$ as in i) respectively ii) generate the
category $\Coh(\kx_\lambda)$, i.e.\ any coherent sheaf on
$\kx_\lambda$ admits a filtration with quotients of this type.
This fact is central for the discussion in \cite{V4}.

Let us now pass from $\Coh$ to $\ka$. Although, $\kx_\lambda$ is
not projective for generic $\lambda$, the category
$\ka_{\kx_\lambda}(\exp(i\omega_\lambda))$ can still be
constructed as in Section \ref{SectAEyieldsDE}.

\begin{proposition}
  For the generic twistor fibre $\kx_\lambda$ the minimal objects
  of the abelian category $\ka_{\kx_\lambda}(\exp(i\omega_\lambda))$ are $k(x)$
  with $x\in
  \kx_\lambda$, and $E[1]$, where $E$ is a $\mu_{\omega_\lambda}$-stable vector
  bundle.

  They generate the abelian category.
\end{proposition}

Note that the second assertion is neither true for $\Coh(X)$ with
$X$ arbitrary (projective or not) nor for $\ka(\exp(i\omega))$ in
the algebraic case. Indeed, the minimal objects $k(x)\in\Coh(X)$
do not generate locally free sheaves neither do the minimal
objects of $\ka(\exp(i\omega))$ described by Proposition
\ref{classmin} generate $F[1]$ with $F$ locally free of negative
slope.

\begin{proof}
For the description of the minimal objects of
$\ka:=\ka_{\kx_\lambda}(\exp(i\omega_\lambda))$ one follows the
arguments in Section \ref{SectMinObj}. The proof that the $k(x)$
and $F[1]$ are minimal did not use any projectivity. To show that
minimal objects are of this form,  observe that all line bundles
are of degree zero. Thus, objects of the form $F[1]$ with $F$
locally free and $\mu$-stable of negative slope do not exist.

Again using that line bundles are of degree zero, one proves that
$\kx_\lambda$ does not contain any curves. Hence, if
$F^\bullet\in\ka$, then $\kh^0(F^\bullet)$ is a torsion sheaf in
dimension zero. For the same reason, $\kh^{-1}(F^\bullet)$ is a
$\mu$-semistable torsion free sheaf of slope zero. Since such a
torsion free sheaf $F$ gives rise to a short exact sequence
$$\xymatrix{0\ar[r]&F\ddual/F\ar[r]&F[1]\ar[r]&F\ddual[1]\ar[r]&0}$$
in $\ka$, every object in $\ka$ does indeed admit a filtration
with quotients isomorphic to $k(x)$ or $E[1]$ as claimed.
\end{proof}

The conclusion of the above discussion is that the minimal
objects of $\ka_X(\exp(i\omega))$ deform naturally to minimal
objects in $\ka_{\kx_\lambda}(\exp(i\omega_\lambda))$ on the
generic twistor fibre, where they generate the abelian category.
In this respect, $\ka_X(\exp(i\omega))$ behaves better than
$\Coh(X)$ itself, as the minimal object, which are simply the
$k(x)$'s, do deform but never generate $\Coh(\kx_\lambda)$.

\medskip

Up to now we have only considered $\ka(\exp(i\omega))$ or, only
virtually more general, $\ka(\exp(B+i\omega))$ with
$(B.\omega)=0$. The case $(B.\omega)$ can be dealt with in a
similar fashion. Roughly, stable vector bundles of slope
$\mu=(B.\omega)$ do not deform sideways in the twistor space, but
their projectivizations do. This is the point of view in
\cite{V4}, which we will complement by briefly discussing the
twistor space associated to the complexified K\"ahler form
$B+i\omega$. This uses the language of Hitchin's generalized
Calabi--Yau structures and their period domains. I believe that
eventually this will be conceptually the right way of dealing
with the general case.

Let us start with a few remarks on the various period domains
(see \cite{HCY} for more details). By definition
$$Q:=\{x~|~(x.x)=0,~(x.\bar
x)>0\}\subset\IP(H^2(X,\IC))\cong\IP^{21}$$ and
$$\tilde Q:=\{x~|~\langle x.x\rangle=0,~\langle x.\bar
x\rangle >0\}\subset\IP(\widetilde H(X,\IC))\cong\IP^{23}.$$
Furthermore, we set $Q':=\tilde Q\cap \IP((H^2\oplus
H^4)(X,\IC))$. These three are the period domains of ordinary K3
surfaces, generalized Calabi--Yau structures $\varphi$, and
generalized Calabi--Yau structures $\varphi$ with $\varphi_0=0$.

For the K3 surface $X$ with its holomorphic volume form $\sigma$
and  with a chosen K\"ahler class $\omega$ one defines
\begin{eqnarray*}
\widetilde T(i\omega)&:=&\IP(\langle{\rm Re}(\sigma),{\rm
Im}(\sigma),{\rm Re}(\exp(i\omega)),{\rm
Im}(\exp(i\omega))\rangle)\cap\tilde Q\\
&\cong&\IP^3\cap\tilde Q.
\end{eqnarray*}
The base of the twistor space $\pi:\kx\to \IP^1$ considered above
is via the period map identified with $T(i\omega):=\widetilde
T(i\omega)\cap Q=\widetilde T(i\omega)\cap Q'$. Note that a
general $\IP^3\subset\IP^{23}$ would intersect $Q$ in only two
points.

For a complexified K\"ahler form $B+i\omega$ one has the
generalized twistor space $\widetilde
T(\exp(i\omega))=\exp(B)\cdot \widetilde T(i\omega)$, which for
$(B.\omega)\ne0$ intersect $Q$ only in two points. The restricted
twistor space $$T(\exp(B+i\omega))= \widetilde
T(\exp(B+i\omega))\cap Q'=\exp(B)\cdot T(i\omega)$$ parametrizes
the generalized Calabi--Yau structures of the form
$\sigma_\lambda+\sigma_\lambda\wedge B$, where $\sigma_\lambda$
is the holomorphic volume form on $\kx_\lambda$.

The above discussion  should in the case $(B.\omega)\ne0$
translate into saying that a $\mu$-stable vector bundle $F$ of
slope $\mu(F)=(B.\omega)$ naturally deforms to a `bundle with
respect to $\sigma_\lambda+\sigma_\lambda\wedge B$'. A general
theory of coherent sheaves for generalized Calabi--Yau structure
(of the form $\sigma+\sigma\wedge B$) still awaits to be
developed, but for rational $B$, they should correspond to
$\alpha_B$-twisted sheaves with $\alpha_B=\exp(B^{0,2})\in
H^2(\ko^*)$. This fits with the point of view in \cite{V4} (see
also \cite[Prop.\ 2.3]{HSch}).

\section{Stability of FM-transforms}

Suppose $\ke$ on $X\times M_H(v)$ is a universal family of
$\mu$-stable vector bundles  on $X$ such that  $M_H(v)$ is
isomorphic to a K3 surface $X'$. The Fourier--Mukai transform with
kernel $\ke$ is an equivalence
$\Phi:=\Phi_\ke:\Db(X)\congpf\Db(X')$. A natural question,
studied in a number of papers (see e.g.\ \cite{B,V,Y}), is the
following:
\begin{itemize}
\item[]{\it When is the image $\Phi(E)$ of a $\mu$-stable vector bundle $E$
again a (shifted) $\mu$-stable vector bundle on $X'$.}
\end{itemize}

(The polarization on $X$ is $H$ and the one on $X'$ has to be
chosen appropriately.) It is known (see \cite{Y}) that one cannot
expect stability in full generality.

The point I wish to make in  this section is that the answer is
yes for the large class of minimal objects in $\ka$ and that our
discussion gives a conceptual and straightforward proof for it.

For the following we let $X'$ be a K3 surface isomorphic to a
fine moduli space $M_H(v)$ of $\mu$-stable vector bundles on $X$
with Mukai vector $v=(r,\ell,s)$. Denote the universal family on
$X\times X'$ by $\ke$ and the induced Fourier--Mukai equivalence
by $\Phi=\Phi_\ke:\Db(X)\congpf\Db(X')$.

\begin{proposition}
  There exists a polarization $H'$ on $X'$ such that for any
  $\mu$-stable vector bundle $E$ on $X$ with $\mu(E)=-(\ell.H)/r$
  one has either
  \begin{itemize}
  \item[$\bullet$]  $\Phi(E)\cong k(y)[-2]$ if $[E\dual]\in M_H(v)$ or otherwise
  \item[$\bullet$] $\Phi(E)\cong F[-1]$ with $F$ a $\mu_{H'}$-stable vector bundle on
  $X'$.
  \end{itemize}
\end{proposition}

\begin{proof} As the proof will show, one could more generally
consider $\mu$-stability with respect to a K\"ahler class $\omega$
(and not $H$) on $X$ and then find a K\"ahler class $\omega'$ on
$X'$ with the asserted property.

First choose a B-field $B\in{\rm NS}(X)_\IQ$ such that the slope
$\mu(\ke_y)=(\ell.\omega)/r$ equals $\beta:=(B.\omega)$. Then
consider the induced stability condition with heart
$\ka_X(\exp(B+i\omega))$.

Following the discussion in Section \ref{SectDA},
$\Phi^H(\exp(B+i\omega))=\lambda\exp(B'+i\omega')$ for some
$\lambda>0$ and a K\"ahler class $\omega'$ on $X'$, which is
rational if $\omega$ was. Furthermore, Proposition \ref{indabequ}
states that $\Psi:=\Phi_{\ke\dual[1]}$ restricts to an equivalence
$$\Psi:\ka_X(\exp(B+i\omega))\congpf\ka_{X'}(\exp(B'+i\omega')).$$
Clearly, under this equivalence minimal objects are mapped to
minimal objects and the minimal objects on either side have been
described by Proposition \ref{classmin}.

If $E$ is a $\mu$-stable vector bundle on $X$ with
$\mu(E)=-\beta$, then $E\dual[1]$ is a minimal object in
$\ka_X(\exp(B+i\omega))$. Thus, either $\Psi(E\dual[1])\cong
k(y)$ for some closed point $y\in X'$ or $\Psi(E\dual[1])\cong
F[1]$ for some $\mu$-stable vector bundle on $X'$. As
$\Psi^{-1}(k(y))=\Phi_{\ke[1]}(k(y))\cong\ke_y[1]$, the first
case occurs precisely if $E\dual\cong\ke_y$ for some $y\in
M_H(v)$, i.e.\ if $[E\dual]\in M_H(v)$.

To conclude, we observe
$\Psi(E\dual[1])\cong\Phi_{\ke\dual[1]}(E\dual[1])\cong
\Phi_{\ke\dual}(E\dual)[2]$ and by Grothendieck--Verdier duality
$\Phi_\ke(E)\cong\Phi_{\ke\dual}(E\dual)\dual[-2]\cong\Psi(E\dual[1])\dual$,
which is either $k(y)\dual\cong k(y)[-2]$ or $(F[1])\dual\cong
F\dual[-1]$.
\end{proof}

So, roughly $\mu$-stable vector bundles of the same slope (up to
sign) as the ones parametrized by the moduli space in question
have $\mu$-stable Fourier--Mukai transform.

\begin{remark}
  i)  In \cite{Y} Yoshioka constructs an explicit example of a
  $\mu$-stable vector bundle with unstable Fourier--Mukai
  transform. An easy check reveals that due to the numerical conditions
  he imposes his example is indeed not covered by the proposition.

  ii) Bartocci et al. have proved in \cite{B} the result for the case
  $\mu=0$. They use the hyperk\"ahler structure of the K3 surface
  and the interpretation of $\mu$-stable vector bundles as bundles on
  the twistor space in a crucial way (see Section \ref{SectTw}).
  In \cite{V} Verbitsky tried
  to generalize this approach to the case of non-vanishing slope
  by working with a similar result for the associated projective
  bundles.

  iii) In \cite{Yos1} Yoshioka presents further results on the
  stability of Fourier--Mukai transforms for the case of Picard
  number one and $\mu(E\dual)=\mu(\ke_y)+1/(\rk(E)r)$. This is not
  covered by our result.
\end{remark}




{\footnotesize \begin{thebibliography}{mm}
\bibitem{B} C.\ Bartocci, U.\ Bruzzo, D.\ Hern\'andez Ruip\'erez
\em A Fourier--Mukai transform for stable bundles on K3 surfaces.
\em J.\ Reine Angew.\ Math.\ 486 (1997), 1-16.

\bibitem{BvdB} A.\ Bondal, M.\ van den Bergh
\em Generators and representability of functors in commutative and
non-commutative geometry. \em Mosc.\ Math.\ J.\ 3 (2003), 1-36.

\bibitem{BSt} T.\ Bridgeland
\em Stability conditions on triangulated categories. \em Ann.\
Math. to appear. math.AG/0212237.

\bibitem{BK3} T.\ Bridgeland
\em Stability conditions on K3 surfaces. \em math.AG/0307164.

\bibitem{Gabriel} P.\ Gabriel
  \em Des cat{\'e}gories ab{\'e}liennes. \em
  Bull.\ Soc.\ Math.\ France  90  (1962), 323-448.

\bibitem{HRS} D.\ Happel, I.\ Reiten, S.\ Smalø
\em Tilting in abelian categories and quasitilted algebras. \em
Memoirs of the AMS 575 (1996).

\bibitem{HLOY} S.\ Hosono, B.H.\ Lian, K.\ Oguiso, S.-T.\ Yau
\em Fourier--Mukai partners of a $K3$ surface of Picard number
one. \em In:  Vector bundles and representation theory (Columbia,
MO, 2002), Contemp.\ Math.\ 322 (2002), 43-55.

\bibitem{HL} D.\ Huybrechts, M.\ Lehn
\em The geometry of moduli spaces of sheaves. \em Aspects of
Math.\ E31. Vieweg, Braunschweig, (1997).


\bibitem{HSch} D.\ Huybrechts, St.\ Schr\"oer
\em The Brauer group of analytic K3 surfaces. \em IMRN.\ 50
(2003), 2687-2698.

\bibitem{HCY} D.\ Huybrechts
\em Generalized Calabi--Yau structures, K3 surfaces, and B-fields.
\em Int.\ J.\ Math.\ 16 (2005), 13-36.

\bibitem{HS} D.\ Huybrechts, P.\ Stellari
\em Equivalences of twisted K3 surfaces. \em Math.\ Ann.\ 332
(2005), 901-936.

\bibitem{HFM} D.\ Huybrechts
\em Fourier--Mukai transforms in algebraic geometry. \em Oxford
Mathematical Monographs (2006).

\bibitem{KS} M.\ Kashiwara, P.\ Schapira
\em Sheaves on manifolds. \em Grundlehren 292. Springer (1990).

\bibitem{Mu} S.\ Mukai
\em On the moduli space of bundles on K3 surfaces, I. \em
 In: Vector Bundles on Algebraic Varieties, Bombay
(1984), 341-413.

\bibitem{OG} K.\ O'Grady
\em The irreducible components of the moduli space of vector
bundles. \em Invent.\ Math.\ 112 (1993), 585-613.

\bibitem{Ogui} K.\ Oguiso
\em K3 surfaces via almost primes. \em Math.\ Res.\ Lett.\ 9
(2002), 47-63.
\bibitem{Or} D.\ Orlov \em On equivalences of derived categories and K3
surfaces. \em J.\ Math.\ Sci.\ (New York) 84 (1997), 1361-1381.

\bibitem{P} A.\ Polishchuk
\em Classification of holomorphic vector bundles on
noncommutative two-tori. \em Documenta Math.\ 9 (2004), 163-181.

\bibitem{St} P.\ Stellari
\em Some remarks about the FM-partners of K3 surfaces with Picard
number 1 and 2. \em Geom.\ Dedicata 108 (2004), 1-13.

\bibitem{Verbitsky}
M.\ Verbitsky \em Hyperholomorphic bundles over a hyper-K\"ahler
manifold. \em J.\ Alg.\  Geom.\  5  (1996),  633-669.

\bibitem{V} M.\ Verbitsky
\em Projective bundles over hyperk\"ahler manifolds and stability
of Fourier--Mukai transform. \em math.AG/0107196.

\bibitem{V4} M.\ Verbitsky
\em Coherent sheaves on general K3 surfaces and tori. \em
math.AG/0205210.

\bibitem{Y3} K.\ Yoshioka
\em Irreducibility of moduli spaces of vector bundles on K3
surfaces. \em math.AG/9907001.


\bibitem{Yos1} K.\ Yoshioka
\em Moduli spaces of stable sheaves on abelian surfaces. \em
Math.\ Ann.\ 321 (2001), 817-884.

\bibitem{Y2} K.\ Yoshioka
\em Twisted stability and Fourier--Mukai transform. \em
math.AG/0106118.

\bibitem{Y} K.\ Yoshioka
\em Stability and the Fourier--Mukai transform. \em Math. Z.\ 245
(2003), 657-665.

\end{thebibliography}}
\end{document}